\documentclass[draft]{amsart}
\usepackage{amsmath}
\usepackage{amssymb}

\newtheorem{thm}{Theorem}
\newtheorem{cor}[thm]{Corollary}
\newtheorem{lem}[thm]{Lemma}

\theoremstyle{remark}
\newtheorem{rmk}{Remark}

\newcommand{\cv}{\mathbb{C}}

\newcommand{\bv}{\mathbb{B}}
\newcommand{\oc}{\mathcal{O}}
\newcommand{\artanh}{\textup{artanh}}
\newcommand{\ds}{\displaystyle}

\begin{document}

\title[The Fridman invariant and the squeezing function]{On the comparison of the Fridman invariant and the squeezing function}

\author[F. Rong, S. Yang]{Feng Rong, Shichao Yang}

\address{School of Mathematical Sciences, Shanghai Jiao Tong University, 800 Dong Chuan Road, Shanghai, 200240, P.R. China}
\email{frong@sjtu.edu.cn}

\address{School of Mathematical Sciences, Shanghai Jiao Tong University, 800 Dong Chuan Road, Shanghai, 200240, P.R. China}
\email{yangshichao@sjtu.edu.cn}

\subjclass[2010]{32H02, 32F45}

\keywords{Fridman invariant, squeezing function, quotient invariant}

\thanks{The authors are partially supported by the National Natural Science Foundation of China (grant no. 11871333).}

\begin{abstract}
Let $D$ be a bounded domain in $\cv^n$, $n\ge 1$. In this paper, we study two biholomorphic invariants on $D$, the Fridman invariant $e_D(z)$ and the squeezing function $s_D(z)$. More specifically, we study the following two questions about the \textit{quotient invariant} $m_D(z)=s_D(z)/e_D(z)$: 1) If $m_D(z_0)=1$ for some $z_0\in D$, is $D$ biholomorphic to the unit ball? 2) Is $m_D(z)$ constantly equal to 1? We answer both questions negatively.
\end{abstract}

\maketitle


\section{Introduction}

Let $D$ be a bounded domain in $\cv^n$, $n\ge 1$, and $\bv^n$ the unit ball in $\cv^n$. Denote by $B_D^k(z,r)$ the $k_D$-ball in $D$ centered at $z\in D$ with radius $r>0$, where $k_D$ is the Kobayashi distance on $D$. Denote by $\oc_u(D_1,D_2)$ the set of \textit{injective} holomorphic maps from $D_1$ into $D_2$.

In \cite{F:invariant}, Fridman introduced the following invariant:
$$h_D(z)=\inf \{1/r:\ B_D^k(z,r)\subset f(\bv^n),\ f\in \oc_u(\bv^n,D),\ f(0)=z\}.$$

For comparison purposes, we call the following variant of $h_D(z)$ the \textit{Fridman invariant} (cf. \cite{DZ, NV:invariants}):
$$e_D(z)=\sup \{\tanh(r):\ B_D^k(z,r)\subset f(\bv^n),\ f\in \oc_u(\bv^n,D),\ f(0)=z\}.$$
Obviously, $e_D(z)=\tanh (h_D^{-1}(z))$ and $h_D(z)=(\artanh(e_D(z)))^{-1}$.

In \cite{DGZ1}, Deng, Guan and Zhang introduced another invariant, called the \textit{squeezing function}, as follows:
$$s_D(z)=\sup \{r:\ r\bv^n\subset f(D),\ f\in \oc_u(D,\bv^n),\ f(z)=0\}.$$

It is clear from the definitions that both the Fridman invariant and the squeezing function are invariant under biholomorphisms, and both take value in $(0,1]$. Both invariants have attracted much attention in recent years, and we refer the readers to the recent survey article \cite{DWZZ:survey} and the references therein for various aspects of the current research on this topic.

Since they are similar in spirit to the Kobayashi-Eisenman volume form $K_D$ and the Carath\'{e}odory volume form $C_D$, respectively, it is a natural problem to consider the comparison of these two biholomorphic invariants. Indeed, recently in \cite{NV:invariants}, Nikolov and Verma have shown the following

\begin{thm}\cite[Proposition 1]{NV:invariants}\label{T:inequality}
Let $D$ be a bounded domain. Then, for any $z\in D$,
$$s_D(z)\le e_D(z).$$
\end{thm}

Regarding the quotient $M_D(z):=C_D(z)/K_D(z)$, the following classical result is well-known (cf. \cite[Theorem E]{W:ball}).

\begin{thm}\cite[Proposition 11.3.3]{K:book}\label{T:quotient}
Let $D$ be a bounded domain in $\cv^n$. If there is a point $z\in D$ such that $M_D(z)=1$, then $D$ is biholomorphic to the unit ball $\bv^n$.
\end{thm}

Obviously, the quotient $m_D(z):=s_D(z)/e_D(z)$ is also a biholomorphic invariant, and Theorem \ref{T:inequality} shows that $m_D(z)\le 1$ for all $z\in D$. It is then natural to ask whether an analogous result to Theorem \ref{T:quotient} holds for $m_D$ (cf. \cite[Concluding Remarks]{MV:invariants}). Our first main result answers this question \textit{negatively}. More precisely, we have the following

\begin{thm}\label{T:main1}
Let $D$ be a bounded, balanced and convex domain in $\cv^n$. Then,
$$m_D(0)=1.$$
\end{thm}

Recall that a domain $D$ is said to be \textit{balanced} if for any $z\in D$, $\lambda z\in D$ for all $|\lambda|\le 1$. And $D$ is said to be \textit{homogeneous} if the automorphism group of $D$ is transitive. As an immediate corollary to Theorem \ref{T:main1}, we have the following

\begin{cor}\label{C:homogeneous}
Let $D$ be a bounded, balanced, convex and homogeneous domain in $\cv^n$. Then,
$$m_D(z)\equiv 1.$$
\end{cor}

On a bounded homogeneous domain $D$, both the Fridman invariant and the squeezing function are constants, which we call the \textit{Fridman constant} $\rho_D^e$ and the \textit{squeezing constant} $\rho_D^s$, respectively. By Corollary \ref{C:homogeneous}, on a bounded, balanced, convex and homogeneous domain $D$, $\rho_D^e=\rho_D^s$, and we simply denote it by $\rho_D$.

It is well-known that the four types of classical Cartan domains are bounded, balanced, convex and homogeneous domains (see e.g. \cite{Mok:book}), and so are their products. Combining Corollary \ref{C:homogeneous} with \cite[Theorem 2]{Ku:Cartan}, we immediately have the following (cf. \cite[Theorem 4.5]{CKK})

\begin{cor}\label{C:product}
Let $D_i$, $1\le i\le m$, be classical Cartan domains and $D=D_1\times D_2\times \cdots \times D_m$. Then,
$$\rho_D=(\rho_{D_1}^{-2}+\rho_{D_2}^{-2}+\cdots+\rho_{D_m}^{-2})^{-1/2}.$$
\end{cor}

Following Theorem \ref{T:inequality} and Corollary \ref{C:homogeneous}, another natural question one can ask is whether $s_D(z)$ is always equal to $e_D(z)$. To our best knowledge, this question is open, mainly because both invariants are hard to compute explicitly. Our second main result answers this question \textit{negatively}. More precisely, we have the following

\begin{thm}\label{T:main2}
On the punctured unit disk $\Delta^\ast$, there exists $\delta>0$ such that
$$m_{\Delta^\ast}(z)<1,\ \ \ \forall\ 0<|z|<\delta.$$
Moreover, $\ds \lim_{z\rightarrow 0} m_{\Delta^\ast}(z)=0$.
\end{thm}

Combining this with the stability results from \cite{F:invariant} and \cite{DGZ2}, we also have the following

\begin{cor}\label{C:annulus}
Let $\delta>0$ be as in Theorem \ref{T:main2}. Fix $z_0\in \Delta^\ast$ with $|z_0|<\delta$. Denote by $A_r$ the annulus $\{z\in \cv:\ r<|z|<1\}$. Then, there exists $\epsilon>0$ such that for all $0<r<\epsilon$, we have
$$m_{A_r}(z_0)<1.$$
\end{cor}

In section \ref{S:equal}, we prove Theorem \ref{T:main1}. In section \ref{S:disk}, we prove Theorem \ref{T:main2} and Corollary \ref{C:annulus}.

\section{Bounded, balanced and convex domains}\label{S:equal}

Let $D$ be a bounded, balanced and convex domain in $\cv^n$. Denote by $k_D$ and $c_D$ the Kobayashi and Carath\'{e}odory distance on $D$, respectively. The following Lempert's theorem is well-known:

\begin{thm}\cite[Theorem 1]{L:convex}\label{T:convex}
On a convex domain $D$, $k_D=c_D$.
\end{thm}

The \textit{Minkowski function} $l_D$ is defined as (see e.g. \cite{JP:book})
$$l_D(z)=\inf \{t>0:\ z/t\in D\},\ \ \ z\in \cv^n.$$
Note that $D=\{z\in \cv^n:\ l_D(z)<1\}$.

Combining Theorem \ref{T:convex} with \cite[Proposition 2.3.1 (c)]{JP:book}, we have the following key lemma.

\begin{lem}\label{L:key}
Let $D$ be a bounded, balanced and convex domain. Then, for all $z\in D$,
$$\tanh (k_D(0,z))=l_D(z).$$
\end{lem}

We also need the following version of the Schwarz Lemma (see e.g. \cite[Theorem 8.1.2]{R:book}).

\begin{lem}\label{L:Schwarz}
Let $D_1$ and $D_2$ be bounded, balanced and convex domains, and $f$ a holomorphic map from $D_1$ to $D_2$ with $f(0)=0$. Then,
$$f(rD_1)\subset rD_2,\ \ \ 0<r\le 1.$$ 
\end{lem}

We now prove Theorem \ref{T:main1} (cf. \cite{A:imbedding}).

\begin{proof}[Proof of Theorem \ref{T:main1}]
Since $s_D(0)\le e_D(0)$ by Theorem \ref{T:inequality}, we only need to show that $e_D(0)\le s_D(0)$.

For any $\epsilon>0$, fix $f\in \oc_u(\bv^n,D)$ with $f(0)=0$, such that $B_D^k(0,r)\subset f(\bv^n)=:\Omega$ and $\tanh(r)>(1-\epsilon)e_D(0)$.

By Lemma \ref{L:key}, we have
$$B_D^k(0,r)=\{z\in D:\ k_D(0,z)<r\}=\{z\in D:\ l_D(z)<\tanh(r)\}.$$
Thus, if we set $L(z):=\tanh(r)z$, then $L$ is a biholomorphism between $D$ and $B_D^k(0,r)$.

Denote by $F$ the holomorphic inverse of $f$ from $\Omega$ to $\bv^n$. If $F(B_D^k(0,r))=\bv^n$, then $D$ is biholomorphic to $\bv^n$ via $F\circ L$, which implies that $s_D(0)=1\ge e_D(0)$.

If $F(B_D^k(0,r))\neq \bv^n$, then let $t\bv^n$ be the maximal ball contained in $F(B_D^k(0,r))$. By the definition of the squeezing function, we have $t\le s_D(0)<1$. Thus, there exists $p\in \bv^n\backslash F(B_D^k(0,r))$ such that $\|p\|<(1+\epsilon)t\le (1+\epsilon)s_D(0)$.

By Lemma \ref{L:Schwarz}, we get
$$(1-\epsilon)e_D(0)<\tanh(r)\le l_D(f(p))\le \|p\|<(1+\epsilon)s_D(0).$$

Since $\epsilon$ is arbitrary, we must have $e_D(0)\le s_D(0)$.
\end{proof}

\begin{rmk}
In the definition of $e_D(z)$, one can also use the Carath\'{e}odory distance $c_D$ instead of the Kobayashi distance $k_D$, and get another related invariant, denoted by $e_D^c(z)$ (cf. \cite{F:polyhedron, NV:invariants}). In \cite[Proposition 1]{NV:invariants}, it was actually proven that $s_D(z)\le e_D^c(z)\le e_D(z)$ for all $z\in D$. When $D$ is a bounded convex domain, by Lempert's theorem, we have $e_D^c(z)=e_D(z)$. Thus, Theorem \ref{T:main1} and Corollary \ref{C:homogeneous} also hold replacing $e_D$ by $e_D^c$.
\end{rmk}

\begin{rmk}
In \cite{Ku:Cartan}, Kubota studied the following invariant:
$$\rho(D):=\sup \{r:\ r\bv^n\subset f(D),\ f\in \oc_u(D,\bv^n)\}.$$
It is clear from the definitions that $\rho(D)=\sup_{z\in D} s_D(z)$. Thus, on a bounded homogeneous domain $D$, we have $\rho(D)=\rho_D^s$.
\end{rmk}

\begin{rmk}
In Corollary \ref{C:product}, if each $D_i$ is the unit disk $\Delta$ then $D=\Delta^n$ is a unit polydisk, and we get
$$\rho_{\Delta^n}=n^{-1/2},\ \ \ h_{\Delta^n}(z)\equiv (\artanh(\rho_D))^{-1}=2\left(\log \frac{\sqrt{n}+1}{\sqrt{n}-1}\right)^{-1},$$
which recovers \cite[Lemma 2.1]{MV:invariants}.
\end{rmk}

\section{The punctured disk and the annulus}\label{S:disk}

First, we prove Theorem \ref{T:main2}.

\begin{proof}[Proof of Theorem \ref{T:main2}]
Let $D=\Delta^\ast$ and $\Omega=\Delta\backslash (-1,0]$. For any fixed $a\in (0,1)$, by the Riemann mapping theorem, there exists a biholomorphic map $f:\Delta\rightarrow \Omega$ with $f(0)=a$. We need to estimate $\tanh(r)$ such that $B_D^k(a,r)\subset f(\Delta)=\Omega$.

By \cite[Corollary 9.1.10]{JP:book}, we have
$$\tanh(k_D(a,z))=\left[\frac{\theta^2+(\log |z|-\log a)^2}{\theta^2+(\log |z|+\log a)^2}\right]^{1/2},\ \ \ \theta=\textup{Arg}(z)\in (-\pi,\pi].$$
In particular, for any $b\in (-1,0)$, setting $x=-\log (-b)$ and $A=-\log a$, we have
$$\tanh(k_D(a,b))=\left[\frac{\pi^2+(x-A)^2}{\pi^2+(x+A)^2}\right]^{1/2}.$$

Consider the function
$$h(x)=\frac{(x-A)^2+\pi^2}{(x+A)^2+\pi^2}.$$
One readily checks that $h(x)$ takes the minimum value at $x=\sqrt{A^2+\pi^2}$. Thus,
$$\begin{aligned}
\min h(x)&=h(\sqrt{A^2+\pi^2})=\frac{(\sqrt{A^2+\pi^2}-A)^2+\pi^2}{(\sqrt{A^2+\pi^2}+A)^2+\pi^2}=1-\frac{4A\sqrt{A^2+\pi^2}}{(\sqrt{A^2+\pi^2}+A)^2+\pi^2}\\
&\ge 1-\frac{4A\sqrt{A^2+\pi^2}}{4A\sqrt{A^2+\pi^2}+\pi^2}=\frac{\pi^2}{4A\sqrt{A^2+\pi^2}+\pi^2}.
\end{aligned}$$
Therefore, there exists $\delta>0$ such that for any $0<a<\delta$, we have
$$e_D(a)^2\ge \frac{\pi^2}{4A\sqrt{A^2+\pi^2}+\pi^2}>e^{-2A}.$$

By \cite[Corollary 7.3]{DGZ1}, we know that $s_D(a)=a=e^{-A}$. Thus, we get $e_D(a)>s_D(a)$ for all $0<a<\delta$. Moreover, the estimate above clearly implies that $\ds \lim_{a\rightarrow 0} m_{\Delta^\ast}(a)=0$.

Since $D$ is circular, this completes the proof.
\end{proof}

\begin{rmk}
In \cite[Lemma 2.2]{MV:invariants}, Mahajan and Verma gave some estimates on $h_{\Delta^\ast}(z)$ as $z\rightarrow 0$. However, we found their estimates hard to use when converted to $e_{\Delta^\ast}(z)$. Thus, we gave a more straightforward estimate above, using \cite[Corollary 9.1.10]{JP:book}.
\end{rmk}

To prove Corollary \ref{C:annulus}, we need the following two stability theorems from \cite{F:invariant} and \cite{DGZ2}, adapted to our setting.

\begin{thm}\cite[Theorem 2.1]{F:invariant}\label{T:stability1}
Let $D$ be a completely hyperbolic domain and $\{D_k\}_{k=1}^\infty$ a sequence of subdomains such that $D_k\subset D_{k+1}$ and $\cup_k D_k=D$. Then, for any $z_0\in D$,
$$\lim_{k\rightarrow \infty} e_{D_k}(z_0)=e_D(z_0).$$
\end{thm}

\begin{thm}\cite[Theorem 2.1]{DGZ2}\label{T:stability2}
Let $D$ be a bounded domain and $\{D_k\}_{k=1}^\infty$ a sequence of subdomains such that $D_k\subset D_{k+1}$ and $\cup_k D_k=D$. Then, for any $z_0\in D$,
$$\lim_{k\rightarrow \infty} s_{D_k}(z_0)=s_D(z_0).$$
\end{thm}

\begin{proof}[Proof of Corollary \ref{C:annulus}]
Suppose that there does not exist $\epsilon>0$ such that for all $0<r<\epsilon$ we have $m_{A_r}(z_0)<1$. Then, there exists a decreasing sequence of $r_k\rightarrow 0$ such that $m_{A_{r_k}}(z_0)=1$ for all $k$. Since $\Delta^\ast$ is completely hyperbolic, $A_{r_k}\subset A_{r_{k+1}}$ and $\cup_k A_{r_k}=\Delta^\ast$, Theorems \ref{T:stability1} and \ref{T:stability2} apply. Thus, $m_{\Delta^\ast}(z_0)=\lim_{k\rightarrow \infty} m_{A_{r_k}}(z_0)=1$, which contradicts Theorem \ref{T:main2}.
\end{proof}

\end{document}